\documentclass[12pt]{article}
\usepackage{amssymb,amsmath}

\input epsf.tex

\title{On non-commutative analytic spaces over non-archimedean fields}

\author {Yan Soibelman}

\begin{document}
\maketitle

\newtheorem{defn}{Definition}
\newtheorem{thm}{Theorem}
\newtheorem{lmm}{Lemma}
\newtheorem{rmk}{Remark}
\newtheorem{prp}{Proposition}
\newtheorem{conj}{Conjecture}
\newtheorem{exa}{Example}
\newtheorem{cor}{Corollary}
\newtheorem{que}{Question}
\newtheorem{ack}{Acknowledgements}
\newcommand{\C}{{\bf C}}
\newcommand{\K}{{\bf k}}
\newcommand{\R}{{\bf R}}
\newcommand{\N}{{\bf N}}
\newcommand{\Z}{{\bf Z}}
\newcommand{\Q}{{\bf Q}}
\newcommand{\G}{\Gamma}
\newcommand{\A}{A_{\infty}}
\newcommand{\g}{{\bf g}}

\newcommand{\epi}{\twoheadrightarrow}
\newcommand{\mono}{\hookrightarrow}
\newcommand\ra{\rightarrow}
\newcommand\uhom{{\underline{Hom}}}
\renewcommand\O{{\cal O}}
\newcommand{\epp}{\varepsilon}

\tableofcontents

\section{Introduction}

\subsection{}
Let $A$ be a unital algebra over the  field of real or complex numbers.
Following [Co] one can think of $A$ as of the algebra of {\it smooth} functions $C^{\infty}(X_{NC})$ on some ``non-commutative
real smooth manifold $X_{NC}$". Differential
geometry on $X_{NC}$ has been developed by Connes and
his followers.
By adding extra structures to $A$ one can define new classes of spaces.
For example if $A$ carries an antilinear involution one can try define a $C^{\ast}$-algebra $C(X_{NC})$
of ``continuous functions on $X_{NC}$ as a completion of $A$ with respect to the norm $|f|=sup_{\pi}||\pi(f)||$, where $\pi$
runs through the set of all topologically irreducible $\ast$-representations in Hilbert spaces.
By analogy with the commutative case, $C(X_{NC})$ corresponds to the non-commutative {\it topological}
space $X_{NC}$. Similarly, von Neumann algebras correspond to non-commutative  measurable
spaces, etc. Main source of new examples for this approach are ``bad" quotients and foliations.

Another class of non-commutative spaces consists of ``non-commutative schemes" and their generalizations.
Here we treat $A$ as the algebra
of {\it regular} functions on the ``non-commutative affine scheme $Spec(A)$". The ground field can be arbitrary (in fact one
can speak about rings, not algebras). Then one can ask whether there is a  ``{\it sheaf}
of regular functions on $Spec(A)$". This leads to the question about localization
in the non-commutative framework.
The localization  of non-commutative associative rings is a
complicated task. An attempt to glue general non-commutative schemes from the affine ones, leads to a replacement of the naively defined category of non-commutative affine schemes by a more
complicated category (see e.g. [R1]).
Main source of examples for this approach is the representation
theory (e.g. theory of quantum groups).

There is an obvious contradiction between two points of view discussed above. Namely, associative  algebras over ${\R}$ or ${\C}$ are treated as algebras of
functions on different types of non-commutative spaces
(smooth manifolds in non-commutative differential geometry and affine schemes
in non-commutative algebraic geometry).
To my knowledge there is no  coherent approach to non-commutative geometry which
resolves this contradiction. In other words,
one cannot start with, say,
non-commutative smooth algebraic variety over ${\C}$, make it
into a non-commutative
complex manifold and then define a non-commutative version of a
smooth structure, so that it becomes a non-commutative
real smooth manifold.
Maybe this is a sign of a general phenomenon: there are many more types
of non-commutative spaces than the commutative ones. Perhaps the traditional
terminology (schemes, manifolds, algebraic spaces, etc.) has to be modified
in the non-commutative world. Although non-commutative spaces resist an attempt
to classify them, it is still
interesting to study
non-commutative analogs of ``conventional" classes of commutative spaces.
Many examples arise if one considers algebras which are
close to commutative (e.g. deformation quantization, quantum groups), or
algebras which are very far from commutative ones (like free algebras).
In a sense these are two ``extreme" cases, and for some reason
the corresponding non-commutative geometry is richer than the ``general" one.

\subsection{}
In this paper we discuss non-commutative analytic spaces
over non-archimedean fields. The list
of ``natural" examples is non-empty (see e.g. [SoVo]).
Analytic non-commutative tori (or elliptic curves) from [SoVo] are different from
$C^{\infty}$ non-commutative tori of Connes and Rieffel.
Although a non-commutative elliptic curve over $\C$ (or over a non-archimedean valuation field, see [SoVo]) appears as a ``bad"
quotient of an analytic space, it carries more
structures than the corresponding ``smooth" bad quotient which is an object of Connes theory.
Non-commutative deformations of a non-archimedean
K3 surface were mentioned in [KoSo1] as  examples of a ``quantization",
which is not formal with respect to the deformation parameter.
It seems plausible that a natural class of ``quantum" non-archimedean
analytic spaces can be derived from cluster ensembles (see [FG]).

Present paper is devoted to  examples
of non-commutative spaces which  can be called {\it non-commutative rigid analytic spaces}.
General theory is far from being developed.
We hope to discuss some of its aspects elsewhere (see [RSo]).

Almost all examples of non-commutative analytic spaces considered in present paper are treated from the point of view of
the approach to rigid analytic geometry offered in
[Be1], [Be2].
The notion of spectrum of
a commutative Banach ring plays a key role in the approach of Berkovich.
Recall that the spectrum $M(A)$ introduced in [Be1] has two equivalent descriptions:

a) the set of multiplicative continuous seminorms on
a unital Banach ring $A$;

b) the set of equivalence classes of
continuous representations of $A$ in
$1$-dimensional Banach spaces over complete Banach fields (i.e. continuous characters).

The space $M(A)$ carries a natural topology so that
it becomes a (non-empty) compact Hausdorff space.
There is a canonical map $\pi: M(A)\to Spec(A)$ which assigns to a multiplicative seminorm
its kernel, a prime ideal. The spectrum $M(A)$ is a natural generalization of
the Gelfand spectrum of a unital commutative $C^{\ast}$-algebra.

Berkovich's definition is very general and does not require $A$ to be an affinoid algebra
(i.e. an admissible quotient of the algebra of analytic functions on a non-archimedean polydisc).
In the affinoid case one can make $M(A)$ into a ringed space (affinoid space).
General analytic spaces are glued from affinoid ones similarly to the gluing of
general schemes from affine schemes.
Classical Tate theory of rigid analytic spaces which is based on the maximal spectrum
of (strictly) affinoid algebras
agrees with Berkovich theory. Analytic spaces in the sense of Berkovich
have better local properties
(e.g. they are locally arcwise connected, see [Be3], while
in the classical rigid analytic geometry the topology is totally disconnected).

Affinoid spaces play the same role of ``building blocks"
in non-archimedean analytic geometry as affine schemes play in the algebraic geometry. For example,
localization of a finite Banach $A$-module $M$
to an affinoid subset $V\subset M(A)$ is achieved by the topological tensoring of $M$ with an affinoid algebra $A_V$, which is the
localization of $A$ on $V$ (this localization {\it is not} an
essentially surjective functor, differently from the case of algebraic geometry).
In order to follow this approach one needs a large supply of ``good" multiplicative subsets of $A$. If $A$ is commutative  this is indeed the case. It is natural to ask what has to be changed in the non-commutative case.

\subsection{}

If $A$ is a non-commutative unital Banach ring (or non-commutative affinoid algebra, whatever this means) then there might be very few ``good" multiplicative
subsets of $A$. Consequently, the supply of affinoid sets can be insufficient to produce a rich theory of non-commutative
analytic spaces. This problem is already known in non-commutative algebraic
geometry, and one can try to look for a possible solution there.
One way to resolve the difficulty was suggested
in [R1]. Namely, instead of localizing {\it rings},
one should localize {\it categories of modules over rings}, e.g. using
the notion of {\it spectrum} of an abelian category introduced in [R1]. Spectrum is a topological space equipped with  Zariski-type topology.
For an associative unital ring $A$ the category of left modules $A-mod$
gives rise to a sheaf of local categories on the spectrum of $A-mod$. If $A$ is commutative,
then the spectrum of $A-mod$ coincides with the usual spectrum $Spec(A)$.
In the commutative case the fiber of the sheaf of categories over $p\in Spec(A)$  is the category of modules over a local ring $A_p$ which is the
localization of $A$ at $p$. In the non-commutative case
the fiber is not a category of modules over a ring.
Nevertheless one can glue general non-commutative spaces from
``affine" ones and call
them {\it non-commutative
schemes}. Thus non-commutative  schemes are  topological spaces equipped with sheaves of local categories (see [R1] for details).

There are more general classes of non-commutative spaces than non-commutative schemes (see e.g. [KoR]). In particular, there might be no underlying topological space (i.e. no ``spectrum"). Then one axiomatizes the notions of covering and descent. Main idea is the following.
If $X=\cup_{i\in I}U_i$ is a ``good" covering of a scheme, then the algebra of functions
$C:={\cal O}(\sqcup_{i,j\in I}(U_i\cap U_j)$ is a coalgebra in the monoidal category
of $A-A$-bimodules, where $A:={\cal O}(\sqcup_{i\in I}U_i)$. In order to have
an equivalence of the category of descent data with the category of quasi-coherent
sheaves on $X$ one deals with the flat topology, which means that $C$ is a (right)
flat $A$-module. In this approach the category of non-commutative spaces
is defined as a localization of the category of coverings with respect to a class
of morphisms called refinements of coverings. This approach can be generalized
to non-commutative case ([KoR]).
One problem mentioned in [KoR] is the absence of a good definition of morphism
of non-commutative spaces defined by means of coverings. The authors developed
another approach based on the same idea, which deals with derived categories of quasi-coherent sheaves rather than with the
abelian categories of quasi-coherent sheaves. Perhaps this approach can be generalized
in the framework of non-commutative analytic spaces discussed in this paper.

\subsection{}

Let me mention some difficulties  one meets
trying to construct a  theory of non-commutative analytic spaces
(some of them are technical but other are conceptual).

1) It is  typical in non-commutative geometry to look
for a point-independent (e.g. categorical) description
of an object or a structure in the commutative case, and then take it as a definition in the
non-commutative case. For example, an affine morphism of schemes can be characterized by the property
that the direct image functor is faithful and exact. This is taken as a definition (see [R1], VII.1.4)
of an affine morphism of non-commutative schemes. Another example is the algebra of regular functions on a quantized
simple group which is defined via Peter-Weyl theorem (i.e. it is defined as the algebra of matrix elements
of finite-dimensional representations of the quantized enveloping algebra, see [KorSo]).
Surprisingly many natural ``categorical" questions do not have
satisfactory (from the non-commutative point of view) answers in analytic case. For example:
how to characterize
categorically an embedding $V\to X=M(A)$, where $V$ is an affinoid domain?

2) In the theory of analytic spaces all rings
are topological (e.g. Banach).  Topology should be involved already in the
definition of the non-commutative version of Berkovich spectrum $M(A)$ as well
as in the localization procedure (question: having a category of coherent sheaves on the Berkovich
spectrum $M(A)$ of an affinoid algebra $A$ how to describe categorically its stalk at a
point $x\in M(A)$?).

3) It is not clear what is a non-commutative analog of the notion
of affinoid algebra.
In the commutative case a typical example of an affinoid algebra is the Tate algebra,
i.e. the completion $T_n$ of the polynomial algebra $K[x_1,...,x_n]$
with respect to the Gauss norm $||\sum_{l\in {\Z}_+}a_lx^l||=sup_l|a_l|$,
where $K$ is a valuation field.
It is important (at least at the level of proofs) that $T_n$ is noetherian.
If we relax the condition that variables $x_i, 1\le i\le n$ commute, then the noetherian
property can fail. This is true, in particular, if one starts
with the polynomial algebra $K\langle x_1,...,x_n\rangle$ of free variables, equipped
with the same Gauss norm as above (it is easy to see that the norm is still
multiplicative).
In ``classical" rigid analytic geometry many proofs are based on the properties of $T_n$
(e.g. all maximal ideals have finite codimension, Weierstrass divison theorem, Noether normalization
theorem, etc.). There is no ``universal"  replacement of $T_n$ in non-commutative world which enjoys the same properties.
On the other hand, there are some candidates which are good for particular
classes of examples. We discuss them in the main body of the paper.

\subsection{}
About the content of the paper. In order to make exposition more transparent
I have borrowed from [RSo] few elementary things, in particular the definition of a non-commutative analog of the Berkovich spectrum $M_{NC}(A)$ (here $A$ is a unital Banach ring). Our definition
is similar to the algebro-geometric definition of
the spectrum $Spec(A):=Spec(A-mod)$ from [R1]. Instead of the category of $A$-modules in [R1]  we  consider here the category
of continuous $A$-modules complete with respect to a seminorm. We prove that $M_{NC}(A)$ is non-empty.
There is a natural map $\pi:M_{NC}(A)\to Spec(A)$. Then we equip $M_{NC}(A)$ with the natural Hausdorff topology.
The set of bounded multiplicative seminorms $M(A)$
is the usual Berkovich spectrum.
Even if $A$ is commutative, the space $M_{NC}(A)$ is larger
than Berkovich spectrum $M(A)$. This phenomenon can be illustrated in the case of non-commutative algebraic
geometry. Instead of considering $k$-points of a commutative ring $A$, where $k$ is a field, one can consider matrix points of $A$, i.e. homomorphisms $A\to Mat_n(k)$. Informally speaking,
such homomorphisms correspond to morphisms of a non-commutative scheme $Spec(Mat_n(k))$ into
$Spec(A)$. Only the case $n=1$ is visible in the ``conventional" algebraic geometry.
Returning to unital Banach algebras we observe that  $M(A)\subset M_{NC}(A)$
regardless of commutativity of $A$. I should say that only $M(A)$ will play
a role in this paper.

Most of the paper is devoted to examples.
We start with the elementary ones (non-commutative polydiscs and
quantum tori). A non-commutative analytic K3 surface
is the most non-trivial example considered in the paper. It can be also called {\it quantum} K3 surface,
because it is a flat deformation  of the analytic K3 surface
constructed in [KoSo1].
We follow the approach of [KoSo1], where
the ``commutative" analytic K3 surface was constructed by gluing it from ``flat"
pieces, each of which is an analytic analog of a Lagrangian
torus fibration in symplectic geometry. General scheme of the construction is explained in Section 7. Main idea is based on the relationship between non-archimedean analytic Calabi-Yau manifolds and real manifolds with integral affine structure discussed in [KoSo1], [KoSo2]. Roughly speaking, to such a Calabi-Yau manifold $X$ we associate
a PL manifold $Sk(X)$, its skeleton. There is a continuos map $\pi:X\to Sk(X)$
such that the generic fiber is an analytic torus. Moreover, there is an embedding
of $Sk(X)$ into $X$, so that $\pi$ becomes a retraction. For the elliptically fibered
K3 surface the skeleton $Sk(X)$ is a $2$-dimensional sphere $S^2$. It has an integral affine
structure which is non-singular outside of the set of $24$ points. It is analogous to the integral
affine structure on the base of a Lagrangian torus fibration in symplectic
geometry (Liouville theorem). Fibers of $\pi$ are Stein spaces. Hence in order
to construct the sheaf ${\cal O}_X$ of analytic functions on $X$ it suffices
to construct $\pi_{\ast}({\cal O}_X)$.

Almost all examples considered in this paper should be called
{\it quantum non-commutative analytic spaces}. They are based on the version
of Tate algebra in which the commutativity of variables $z_iz_j=z_jz_i$
is replaced by the $q$-commutativity $z_iz_j=qz_jz_i, i<j, q\in K^{\times}, |q|=1$.
In particular our non-commutative analytic K3 surface is defined as a ringed space, with the
underlying topological space being an ordinary K3 surface equipped with the natural
Grothendieck topology introduced in [Be3] and the sheaf of non-commutative algebras
which is locally isomorphic to a quotient of the above-mentioned
``quantum" Tate algebra.
The construction of a quantum K3 surface
uses a non-commutative analog of the map $\pi$. In othere words, the skeleton
$Sk(X)$ survives under the deformation procedure.
We will point out a more general phenomenon, which does not
exist in ``formal" deformation quantization.
Namely, as the deformation parameter $q$ gets closer to $1$
we recover more and more of the Berkovich spectrum
of the undeformed algebra.

\subsection{}

The theory of non-commutative non-archimedean analytic spaces
should have applications to mirror symmetry in the spirit of [KoSo1],[KoSo2].
More precisely, it looks plausible that certain deformation of the Fukaya category of a
maximally degenerate (see [KoSo1], [KoSo2]) hyperkahler manifold can be
realized as the derived category of coherent sheaves on the non-commutative deformation
of the dual Calabi-Yau manifold (which is basically the same hyperkahler manifold).
This remark was one of the motivations of this paper. Another motivation is the theory of $p$-adic
quantum groups. We will discuss it elsewhere.

{\it Acknowledegements}. I am grateful to IHES and Max-Planck Institut f\"ur Mathematik
for excellent research conditions. I thank to Vladimir Berkovich,
Vladimir Drinfeld, Maxim Kontsevich and
Alexander Rosenberg for useful discussions. I thank to Hanfeng Li for comments on the
paper.
I am especially grateful to Lapchik (a.k.a Lapkin) for encouragement, constant interest and multiple remarks on the preliminary drafts (see [Lap]). Main object of present paper should be called ``lapkin spaces".
Only my unability to discover the  property ``sweetness" (predicted by Lapchik,
cf. with the ``charm" of quarks) led me to the (temporary) change of terminology to
``non-commutative non-archimedean analytic spaces".

\section{Non-commutative Berkovich spectrum}

\subsection{Preliminaries}
We refer to [Be1], Chapter 1 for the terminology of seminormed groups, etc. Here we recall few terms for convenience of the reader.

Let $A$ be an associative unital Banach ring. Then, by definition, $A$ carries a norm
$|\bullet|_A$, and moreover, $A$ is
complete with respect to this norm. The norm is assumed to be submultiplicative, i.e. $|ab|_A\le |a|_A|b|_A,
a,b\in A$ and unital, i.e. $|1|_A=1$. We call the norm {\it non-archimedean} if,
instead of the usual inequality $|a+b|_A\le |a|_A+|b|_A, a,b\in A$, we have a stronger one
$|a+b|_A\le max\{|a|_A,|b|_A\}$. A {\it seminormed module} over $A$ is (cf. with [Be1])
a left unital (i.e. $1$ acts as $id_M$) $A$-module $M$ which carries a seminorm $|\bullet|$ such that $|am|\le C|a|_A|m|$ for
some $C>0$ and all $a\in A,m\in M$. Seminormed $A$-modules form a category $A-mod^c$, such that
a morphism $f:M\to N$ is a homomorphism of $A$-modules satisfying the condition
$|f(m)|\le const\,|m|$ (i.e. $f$ is bounded). Clearly the kernel $Ker(f)$ is closed
with respect to the topology defined by the seminorm on $M$.
A morphism $f:M\to N$ is called {\it admissible}
if the quotient seminorm on $Im(f)\simeq M/Ker(f)$ is equivalent to the one induced from $N$.
We remark that the category $A-mod^c$ is not abelian in general.
Following [Be1] we call {\it valuation field} a commutative Banach field $K$ whose norm
is multiplicative, i.e. $|ab|=|a||b|$. If the norm is non-archimedean, we will cal $K$ a {\it non-archimedean} valuation field. In this case
one can introduce a valuation map $val:K^{\times}\to \R\cup {+\infty}$,
such that $val(x)=-log|x|$.

\subsection{Spectrum of a Banach ring}

Let us introduce a partial order on the objects of $A-mod^c$. We say that
$N\ge_cM$ if there exists a closed admissible embedding $i:L\to \oplus_{I-finite}N$
and an admissible epimorphism $pr:L\to M$. We will denote by $N\ge M$ a similar
partial order on the category $A-mod$ of left $A$-modules (no admissibility condition
is imposed). We will denote by $=_c$ (resp. $=$ for $A-mod$)
the equivalence relation generated by the above partial order.

Let us call an object of $A-mod^c$ (resp. $A-mod$) {\it minimal} if it satisfies
the following conditions (cf. [R1]):

1) if $i:N\to M$ is a closed admissible embedding (resp. any embedding in the case of $A-mod$)
then $N\ge_cM$ (resp. $N\ge M$);

2) there is an element $m\in M$ such that $|m|\ne 0$ (resp. $m\ne 0$ for $A-mod$).

We recall (see [R1]) that the spectrum $Spec(A):=Spec(A-mod)$ consists of equivalence classes (w.r.t. $=$) of minimal objects
of $A-mod$. It is known that $Spec(A)$ is non-empty and contains classes of $A$-modules
$A/m$, where $m$ is a left maximal ideal of $A$. Moreover, $Spec(A)$ can be identified with
the so-called {\it left spectrum} $Spec_l(A)$, which can be also described in terms of
a certain subset of the set of
left ideals of $A$ (for commutative $A$ it is the whole set of prime ideals).

Let $A-mod^b$ be a full subcategory of $A-mod^c$ consisting of $A$-modules which are complete
with respect to their seminorms (we call them {\it Banach modules} for short).

\begin{defn} The non-commutative analytic spectrum $M_{NC}(A)$ consists of classes of equivalence (w.r.t. to $=_c$)
of minimal (i.e satisfying 1) and 2))
objects $M$ of $A-mod^b$ which satisfy also the following property:

3) if $m_0\in M$ is such that $|m_0|\ne 0$ then  the left module $Am_0$ is minimal, i.e.
defines a point of $Spec(A)$ (equivalently, this means that $p:=Ann(m_0)\in Spec_l(A)$),
and this point does not depend on a choice of $m_0$.

\end{defn}

The following easy fact implies that $M_{NC}(A)\ne \emptyset$.

\begin{prp} Every (proper) left maximal ideal $m\subset A$ is closed.

\end{prp}

{\it Proof.} We want to prove that the  closure $\overline{m}$ coincides with $m$.
The ideal $\overline{m}$ contains $m$. Since $m$ is maximal, then either
$\overline{m}=m$ or $\overline{m}=A$. Assume the latter.
Then there exists a sequence $x_n\to 1, n\to +\infty, x_n\in m$. Choose $n$ so large that
$|1-x_n|_A<1/2$. Then $x_n=1+(x_n-1):=1+y_n$ is invertible in $A$,
since $(1+y_n)^{-1}=\sum_{l\ge 0}(-1)^ly_n^l$ converges. Hence ${m}=A$.
Contradiction. $\blacksquare$

\begin{cor} $M_{NC}(A)\ne \emptyset$

\end{cor}

{\it Proof.} Let $m$ be a left maximal ideal in $A$ (it does exists because of the standard arguments
which use Zorn lemma). It is closed by previous Proposition.
Then $M:=A/m$ is a cyclic Banach $A$-module with respect to the quotient norm.
We claim that it contains no proper closed submodules.
Indeed, let $N\subset M$ be a proper closed submodule. We may assume it contains an element
$n_0$ such that annihilator $Ann(n_0)$ contains $m$ as a proper subset. Since $m$ is maximal
we conclude that $Ann(n_0)=A$. But this cannot be true since $1\in A$ acts without kernel
on $A/m$, hence $1\notin Ann(n_0)$. This contradiction shows that
$M$ contains no proper closed submodules. In order to finish the proof we recall that
simple $A$-module $A/m$ defines a point of $Spec(A)$. Hence the conditions 1)-3) above are satisfied. $\blacksquare$

Abusing terminology we will often say that an object belongs to the spectrum (rather than saying
that its equivalence class belongs to the spectrum).

\begin{prp} If $M\in A-mod^b$ belongs to $M_{NC}(A)$ then its seminorm is, in fact,
a Banach norm.

\end{prp}

{\it Proof.} Let $N=\{m\in M||m|=0\}$. Clearly $N$ is a closed submodule. Let $L$
be a closed submodule of the finite sum of copies of $N$, such that there exists
an admissible epimorphism $pr:L\to M$.
Then the submodule $L$ must carry trivial induced
seminorm, and, moreover, admissibility of the epimorphism $pr:L\to M$ implies that
the seminorm on $M$ is trivial. This contradicts to 2). Hence $N=0$. $\blacksquare$

\begin{rmk} a) By condition 3) we have a map of sets $\pi:M_{NC}(A)\to Spec(A)$.

b) If $A$ is commutative then any bounded multiplicative seminorm on $A$ defines
a prime ideal $p$ (the kernel of seminorm). Moreover, $A/p$ is a Banach $A$-module,
which belongs to $M_{NC}(A)$. Hence Berkovich spectrum $M(A)$ (see [Be1])
is a subset of $M_{NC}(A)$. Thus $M_{NC}(A)$ can be also called non-commutative
Berkovich spectrum.

\end{rmk}

\subsection{Topology on $M_{NC}(A)$}

If a Banach $A$-module $M$ belongs to $M_{NC}(A)$
then it is equivalent (with respect to $=_c$) to the closure of any cyclic submodule
$M_0=Am_0$. Then for a fixed $a\in A$
we have a function $\phi_a:(M,|\bullet|)\mapsto |am_0|$, which
can be thought of as a real-valued function on $M_{NC}(A)$.

\begin{defn} The topology on $M_{NC}(A)$ is taken to be the weakest one for which
all functions $\phi_a, a\in A$ are continuous.

\end{defn}

\begin{prp} The above topology makes $M_{NC}(A)$ into a Hausdorff topological space.

\end{prp}

{\it Proof.} Let us take two different points of the non-commutative
analytic spectrum $M_{NC}(A)$ defined by cyclic seminormed modules $(M_0, |\bullet|)$
and $(M_0^{\prime}, |\bullet|^{\prime})$. If $M_0$ is not equivalent to $M_0^{\prime}$ with respect
to $=$ (i.e. in $A-mod$),
then there are exist two different closed left ideals $p,p^{\prime}$ such that
$M_0\simeq A/p$ and $M_0^{\prime}\simeq A/p^{\prime}$ (again,
this is an isomorphism of $A$-modules only. Banach norms
are not induced from $A$). Then there is an element $a\in A$ which belongs, to,say, $p$ and
does not belong to $p^{\prime}$ (if $p\subset p^{\prime}$ we interchange $p$ and $p^{\prime}$).
Hence the function $\phi_a$ is equal to zero on the closure of $(M_0, |\bullet|)$ and
$\phi_a(\overline{(M_0^{\prime},|\bullet|^{\prime}})=c_a>0$. Then open sets
$U_0=\{x\in M_{NC}(A)| \phi_a(x)<c_a/2\}$ and
$U_0^{\prime}=\{x\in M_{NC}(A)| 3c_a/4<\phi_a(x)<c_a\}$ do not intersect and separate two given
points of the analytic spectrum.

Suppose $M=M_0=M_0^{\prime}=Am_0$.  Since the corresponding points of $Spec(A)$ coincide,
we have the same cyclic $A$-module which carries two different norms $|\bullet|$ and
$|\bullet|^{\prime}$. Let  $am_0=m\in M$ be an element such that $|m|\ne |m|^{\prime}$.
Then the function $\phi_a$ takes different values at the corresponding
points of $M_{NC}(A)$ (which are completions of $M$ with respect to the above norms),
and we can define separating open subsets as before. This concludes
the proof. $\blacksquare$

\subsection{Relation to multiplicative seminorms}

If $x:A\to {\R}_+$ is a multiplicative bounded seminorm on $A$ (bounded means $x(a)\le |a|_A$
for all $a\in A$) then $Ker\,x$ is a closed $2$-sided ideal in $A$. If $A$ is non-commutative, it can contain
very few $2$-sided ideals. At the same time, a bounded multiplicative
seminorm $x$  gives rise to a point of $M_{NC}(A)$ such that
the corresponding Banach $A$-module is $A/Ker\,x$ equipped with the left action of $A$. 

There exists a class of {\it submultiplicative} bounded seminorms
on $A$ which is contained in $M(A)$, if $A$ is commutative. More precisely,
let us consider the set $P(A)$  of all submultiplicative bounded seminorms on $A$.
By definition, an element of $P(A)$ is a seminorm such that
$|ab|\le |a||b|, |1|=1, |a|\le C |a|_A$ for all $a,b\in A$ and some $C>0$ ($C$ depends
on the seminorm). The set $P(A)$ carries natural partial order:
$|\bullet|_1\le |\bullet|_2$ if $|a|_1\le |a|_2$ for all $a\in A$.
Let us call {\it minimal seminorm} a minimal element of $P(A)$ with respect to
this partial order, and denote by $P_{min}(A)$ the subset of minimal seminorms.
The latter set is non-empty by Zorn lemma. Let us recall the following
classical result (see e.g. [Be1]).

\begin{prp} If $A$ is commutative then $P_{min}(A)\subset M(A)$, i.e. any minimal
seminorm is multiplicative.

\end{prp}

Since there exist multiplicative bounded seminorms which are not minimal
(take e.g. the trivial seminorm on the ring of integers ${\Z}$ equipped with the usual absolute value Banach
norm) it is not reasonable to define $M(A)$ in the non-commutative case
as a set of minimal bounded seminorms. On the other hand one can prove
the following result.

\begin{prp}
If $A$ is left noetherian as a ring then a minimal bounded seminorm defines
a point of $M_{NC}(A)$.

\end{prp}

{\it Proof.} Let $p$ be the kernel of a minimal seminorm $v$.
Then $p$ is a $2$-sided closed ideal. The quotient $B=A/p$ is a Banach
algebra with respect to the induced norm. It is topologically simple,
i.e. does not contain non-trivial closed $2$-sided ideals. Indeed,
let $r$ be a such an ideal. Then $B/r$ is a Banach algebra with
the norm induced from $B$. The pullback of this norm to $A$
gives rise to a bounded seminorm on $A$, which is
smaller than $v$, since it is equal to zero on a closed $2$-sided
ideal which contains $p=Ker\,v$. The remaining proof
that  $A/p\in M_{NC}(A)$ is similar to the one from [R1]. Recall that it
was proven in [R1] that if $n$ is a $2$-sided ideal in
a noetherian ring $R$ such that $R/n$ contains no $2$-sided ideals
then $A/n\in Spec(A-mod)$. $\blacksquare$

We will denote by $M(A)$ the subset of $M_{NC}(A)$ consisting
of bounded multiplicative seminorms. It carries the induced topology,
which coincides for a commutative $A$ with the topology introduced in [Be1]. We will call the corresponding topological space
{\it Berkovich spectrum} of $A$.

Returning to the beginning of this section we observe that the set $P(A)$
of submultiplicative bounded seminorms contains $M_{NC}(A)$.
Indeed if $v\in P(A)$ then we have a left Banach $A$-module $M_v$,
which is the completion of $A/Ker\,v$ with respect to the
norm induced by $v$. Thus $M_v$ is a cyclic Banach $A$-module.
A submultiplicative bounded seminorm which defines a
point of the analytic spectrum can be characterized by the following property:
$v\in P(A)$ belongs to $M_{NC}(A)$ iff $A/Ker\,v \in Spec(A-mod)$
(equivalently, if $Ker\,v\in Spec_l(A)$).

\subsection{Remark on representations in a Banach vector space}

Berkovich spectrum of a commutative unital  Banach ring $A$ can be understood
as a set of equivalence classes of one-dimensional representations
of $A$ in Banach vector spaces over valuation fields. 
One can try to do a similar thing in the non-commutative
case based on the following simple considerations.

Let $A$ be as before, $Z\subset A$ its center. Clearly it is a commutative unital Banach subring
of $A$. It follows from [Be1], 1.2.5(ii) that there exists a bounded seminorm on $A$ such that
its restriction to $Z$ is multiplicative (i.e. $|ab|=|a||b|$). Any such a seminorm $x\in M(Z)$
gives rise to a valuation field $Z_x$, which is the completion (with respect to the induced
multiplicative norm) of the quotient field of the domain $Z/Ker\,x$. Then the completed tensor
product $A_x:=Z_x\widehat{\otimes}_ZA$ is a Banach $Z_x$-algebra.
For any valuation field $F$ and a Banach $F$-vector space $V$ we will denote by
$B_F(V)$ the Banach algebra of all bounded operators on $V$.
Clearly the left action of $A_x$ on itself is continuous. Thus we have a homomorphism
of Banach algebras $A_x\to B_{Z_x}(A_x)$. Combining this homomorphism with the homomorphism
$A\to A_x, a\mapsto 1\otimes a$ we see that the following result holds.

\begin{prp} For any associative unital Banach ring $A$ there exists a
valuation field $F$, a Banach $F$-vector space $V$ and a representation
$A\to B_F(V)$ of $A$ in the algebra of bounded linear operators in $V$.

\end{prp}

\section{Non-commutative affinoid algebras}

Let $K$ be a non-archimedean valuation field, $r=(r_1,...,r_n)$,
where $r_i>0,1\le i\le n$.
In the ``commutative" analytic geometry an affinoid algebra $A$ is defined
as an admissible quotient of the unital Banach algebra $K\{T\}_r$ of formal series
$\sum_{l\in {\Z}_+^n}a_lT^l$, such that $max\,|a_l|r^l\to 0,l=(l_1,...,l_n)$. The latter algebra
is the completion of the algebra of polynomials $K[T]:=K[T_1,...,T_n]$
with respect to
the norm $|\sum_la_lT^l|_r=max|a_l|r^l$. In the non-commutative
case we start with the algebra
$K\langle T\rangle:=K\langle T_1,...,T_n\rangle$ of polynomials
in $n$ free variables
and consider its completion $K\langle \langle T\rangle\rangle_r$ with respect to the norm
$|\sum_{\lambda\in P({\Z}_+^n)}a_{\lambda}T^{\lambda}|=max_{\lambda} |a_{\lambda}|r^{\lambda}$.
Here $P({\Z}_+^n)$ is the set of finite paths in ${\Z}_+^n$ starting
at the origin, and $T^{\lambda}=T_1^{\lambda_1}T_2^{\lambda_2}....$ for the path
which moves $\lambda_1$ steps in the direction $(1,0,0,...)$ then $\lambda_2$
steps in the direction $(0,1,0,0...)$, and so on (repetitions are allowed, so we
can have a monomial like $T_1^{\lambda_1}T_2^{\lambda_2}T_1^{\lambda_3}$).
We say that a Banach unital algebra $A$ is {\it non-commutative affinoid algebra}
if there is an admissible surjective homomorphism
$K\langle \langle T\rangle\rangle_r\to A$. In particular, $K\{T\}_r$ is such an algebra,
and hence all commutative affinoid algebras are such algebras.
We will restrict ourselves to the class of {\it noetherian} non-commutative
affinoid algebras, i.e. those which are noetherian as left rings. All classical
affinoid algebras belong to this class.

For noetherian affinoid algebras one can prove the following result (the proof is similar
to the proof of Prop. 2.1.9 from [Be1], see also [SchT]).

\begin{prp} Let $A$ be a noetherian non-commutative affinoid algebra.
Then the category $A-mod^f$ of (left) finite $A$-modules is equivalent
to the category $A-mod^{fb}$ of (left) finite Banach $A$-modules.

\end{prp}

An important class of noetherian affinoid algebras consists of
{\it quantum} affinoid algebras. By definition they are
admissible quotients of algebras $K\{T_1,...,T_n\}_{q,r}$. The latter
consists of formal power series
$f=\sum_{l\in {\Z}_+^n}a_lT^l$ of $q$-commuting variables
(i.e. $T_iT_j=qT_jT_i, i<j$ and $q\in K^{\times}, |q|=1$) such that
$|a_l|r^l\to 0.$ Here $T^l=T_1^{l_1}...T_n^{l_n}$ (the order
is important  now).
It is easy to see that for any
$r=(r_1,...,r_n)$ such that all $r_i> 0$ the function $f\mapsto |f|_r:=max_l|a_l|r^l$ defines
a multiplicative norm on the polynomial algebra $K[T_1,...,T_n]_{q,r}$
in $q$-commuting variables $T_i, 1\le i\le n$. Banach algebra
$K\{T_1,...,T_n\}_{q,r}$ is the completion of the latter with respect to the norm $f\mapsto |f|_r$.
Similarly, let $Q=((q_{ij}))$ be an $n\times n$ matrix with
entries from $K$ such that $q_{ij}q_{ji}=1, |q_{ij}|=1$ for all $i,j$.
Then we define the quantum affinoid algebra
$K\{T_1,...,T_n\}_{Q,r}$ in the same way as above, starting
with polynomials in variables $T_i, 1\le i\le n$
such that $T_iT_j=q_{ij}T_jT_i$.
One can  think of this Banach algebra as of the
quotient of
$K\langle\langle T_i,t_{ij}\rangle \rangle_{r,{\bf 1}_{ij}}$,
where $1\le i,j\le n$ and ${\bf 1}_{ij}$ is the unit $n\times n$
matrix,
by the  two-sided ideal generated by the relations
$$t_{ij}t_{ji}=1, T_iT_j=t_{ij}T_jT_i, t_{ij}a=at_{ij},$$
for all indices $i,j$ and all
$a\in K\langle\langle T_i,t_{ij}\rangle \rangle_{r,{\bf 1}_{ij}}$.
In other words, we treat $q_{ij}$ as variables
which belong to the center of our algebra and have the
norms equal to $1$.

\section{Non-commutative analytic affine spaces}

Let $k$ be a commutative unital Banach ring.
Similarly to the previous section
we start with the algebra ${k}\langle\langle x_1,...,x_n\rangle\rangle$
of formal series in free variables $T_1,...,T_n$.
Then for each $r=(r_1,...,r_n)$
we define a subspace ${k}\langle\langle T_1,...,T_n\rangle\rangle_r$
consisting of series $f=\sum_{i_1,...,i_m}a_{i_1,...,i_m}T_{i_1}...T_{i_m}$
such that $\sum_{i_1,...,i_m}|a_{i_1,...,i_m}|r_{i_1}...r_{i_m}<+\infty$.
Here the summation is taken over all sequences $(i_1,...,i_m), m\ge 0$ and
$|\bullet|$ denotes the norm in $k$.
In the non-archimedean case the convergency condition is replaced
by the following one: $max\,|a_{i_1,...,i_m}|r_{i_1}...r_{i_m}<+\infty$.
Clearly each ${k}\langle\langle x_1,...,x_n\rangle\rangle_r$ is a Banach
algebra (called the algebra of analytic functions on a non-commutative
$k$-polydisc $E_{NC}(0,r)$, cf. [Be1], 1.5). We
would like to define a non-commutative
$n$-dimensional analytic $k$-affine space ${\bf A}^n_{NC}$ as the coproduct
$\cup_rE_{NC}(0,r)$. By definition the algebra of analytic functions on the quantum affine space
is given by the above series such that
$max\,|a_{i_1,...,i_m}|r_{i_1}...r_{i_m}<+\infty$ for all $r=(r_1,...,r_n)$.
In other words, analytic functions are given by the series which are convergent
in all non-commutative polydics with centers in the origin.

The algebra of analytic functions on the non-archimedean {\it quantum closed
polydisc} $E_{q}(0,r)$ is, by definition, $k\{T_1,...,T_n\}_{q,r}$.
The algebra of analytic functions on non-archimedean {\it quantum affine space ${\bf A}^N_q$}
consists of the  series $f$ in $q$-commuting variables $T_1,...,T_n$, such that for all $r$ the norm
$|f|_{r}$ is finite. Equivalently, it is the coproduct of {\it quantum closed
polydiscs} $E_{q}(0,r)$.
There is an obvious generalization of this example
to the case when $q$ is replaced by a matrix $Q$ as in
the previous section.
We will keep the terminology for the matrix case as well.

\section{Quantum analytic tori}

Let $K$ be a non-archimedean valuation field, $L$ a free abelian
group of finite rank $d$, $\varphi:L\times L\to {\Z}$ is a skew-symmetric
bilinear form, $q\in K^{\ast}$ satisfies the
condition $|q|=1$. Then $|q^{\varphi(\lambda,\mu)}|=1$
for any $\lambda,\mu \in L$. We denote by
$A_q(T(L,\varphi))$ the
{\it algebra of regular functions on the quantum torus $T_q(L,\varphi)$}.
By definition, it is a $K$-algebra with generators $e(\lambda),\lambda\in L$,
subject to the relation
$$e(\lambda)e(\mu)=q^{ \varphi(\lambda,\mu)}e(\lambda+\mu).$$

\begin{defn} The space ${\cal O}_q(T(L,\varphi,(1,...,1)))$ of analytic functions
on the quantum torus of multiradius $(1,1,...,1)\in \Z_+^d$  consists of series
$\sum_{\lambda\in L}a(\lambda)e(\lambda)$,$a(\lambda)\in K$ such that
$|a(\lambda)|\to 0$ as $|\lambda|\to \infty$ (here $|(\lambda_1,...,\lambda_d)|=
\sum_i|\lambda_i|$).

\end{defn}

It is easy to check (see [SoVo]) the following result.

\begin{lmm} Analytic functions ${\cal O}_q(T(L,\varphi,(1,...,1)))$ form
a Banach $K$-algebra. Moreover, it is a  noetherian quantum
affinoid algebra.
\end{lmm}

This example admits the following generalization.
Let us fix a basis $e_1,...,e_n$ of $L$ and positive
numbers $r_1,...,r_n$. We define
$r^{\lambda}=r_1^{\lambda_1}...r_n^{\lambda_n}$ for
any $\lambda=\sum_{1\le i\le n}\lambda_ie_i\in L$.
Then the algebra ${\cal O}_q(T(L,\varphi,r))$ of analytic functions on the
{\it quantum torus of multiradius $r=(r_1,...,r_n)$} is defined
by same series as in Definition 3, with the only change that
$|a(\lambda)|r^{\lambda}\to 0$ as $|\lambda|\to \infty$.
We are going to denote the corresponding non-commutative analytic
space by $T^{an}_q(L,\varphi,r)$. It is a quantum affinoid space.
The coproduct $T^{an}_q(L,\varphi)=\cup_r T^{an}_q(L,\varphi,r)$ is called the
{\it quantum analytic torus}. The algebra of analytic functions
${\cal O}_q(T(L,\varphi))$
on it consists, by definition, of the above series such
that for all $r=(r_1,...,r_d)$ one has:
$|a(\lambda)|r^{\lambda}\to 0$ as $|\lambda|\to \infty.$
To be consistent with the notation of the previous subsection we will
often denote the dual to $e_i$ by $T_i$.

\subsection{Berkovich spectra of quantum polydisc and analytic torus}

Assume for simplicity that the pair $(L,\varphi)$ defines a
simply-laced lattice of rank $d$, i.e. for some basis $(e_i)_{1\le i\le d}$
of $L$ one has $\varphi(e_i,e_j)=1, i<j$. In the coordinate notation we have
$T_iT_j=qT_jT_i, i<j, |q|=1$.
Any $r=(r_1,...,r_d)\in \R_{>0}^d$ gives rise to a point
$\nu_r\in M({\cal O}_q(T(L,\varphi)))$ such that $\nu_r(f)=
max_{\lambda \in L}|a(\lambda)|r^{\lambda}$. In this way we obtain a (continuous)
embedding of $\R_{>0}^d$ into Berkovich spectrum of quantum torus.
This is an example of a more general phenomenon.
In fact $\R_{>0}^d$ can be identified with
the so-called {\it skeleton} $S((G_m^{an})^d)$
(see [Be1], [Be3])
of the $d$-dimensional (commutative) analytic torus $(G_m^{an})^d$.
Skeleta can be defined
for more general analytic spaces (see [Be3]). For example, the skeleton of the $d$-dimensional
Drinfeld upper-half space $\Omega^d_K$ is the Bruhat-Tits building of $PGL(d,K)$. There is also a different
notion of skeleton,
which makes sense for so-called maximally degenerate
Calabi-Yau manifolds (see [KoSo1], Section 6.6 for the details).
In any case a skeleton is a PL-space (polytope) naturally
equipped with the sheaf of affine functions.
One can expect that the notion of skeleton
(either in the sense of Berkovich or in the sense
of [KoSo1]) admits a generalization
to the case of quantum analytic spaces modeled by quantum affinoid algebras.

We have constructed above
an embedding of  $S((G_m^{an})^d)$ into
$M({\cal O}_q(T(L,\varphi))$.
If $q=1$ then there is a retraction
$(G_m^{an})^d\to S((G_m^{an})^d)$.
The pair $((G_m^{an})^d, S((G_m^{an})^d))$
is an example of {\it analytic torus fibration} which plays an important role in
mirror symmetry (see [KoSo1]). One can expect
that this picture
admits a quantum analog.

The skeleton of the analytic  quantum torus survives
in $q$-deformations as long as $|q|=1$. Another example of
this sort is a quantum K3 surface considered in Section 7.
One can expect that the skeleton survives under $q$-deformations
with $|q|=1$
for all analytic spaces
which have a skeleton. Then it is natural to ask
whether the
Berkovich spectrum of a quantum non-archimedean analytic space
contains more than just the skeleton.
Surprisingly, as $q$ gets sufficiently close to $1$, the answer is positive.
Let  $\rho=(\rho_1,...,\rho_d),r=(r_1,...,r_d)\in \R_{\ge 0}^d$ and $a=(a_1,...,a_d)\in K^d$. We
assume that $|1-q|<1$ and $|a|\le |\rho|<r$ (i.e. $a_i\le \rho_i,<r_i, 1\le i\le d$).
Let $f=\sum_{n\in \Z_+^d}c_nT^n$ be a polynomial in $q$-commuting variables. Set $t_i=T_i-a_i, 1\le i\le d$.
Then $f$ can be written as $f=\sum_{n\in \Z_+^d}b_nt^n$ (although $t_i$ are no longer $q$-commute).

\begin{prp} The seminorm
$\nu_{a,\rho}(f):=max_n |b_n|\rho^n$
defines a point of the Berkovich spectrum of the quantum
closed polydisc disc $E_q(0,r)$.

\end{prp}

{\it Proof.} It suffices to show that $\nu_{a,r}(fg)=\nu_{a,r}(f)\nu_{a,r}(g)$ for any
two polynomials $f,g\in K[T_1,...,T_d]_q$.

Let us introduce new variables $t_i$ by the formulas $T_i-a_i=t_i, 1\le i\le d$. Clearly $\nu_{a,\rho}(t_i)=\rho_i, 1\le i\le d$.
By definition $t_it_j=qt_jt_i+(q-1)(a_it_j+a_jt_i)+(q-1)a_ia_j$.
Therefore, for any multi-indices $\alpha,\beta$ one
has $t^{\alpha}t^{\beta}=q^{\varphi(\alpha,\beta)}t^{\beta}t^{\alpha}+D$,
where $\nu_{a,\rho}(D)<\rho^{|\alpha|+|\beta|}$. Here
$t^{(\alpha_1,...,\alpha_d)}:=t_1^{\alpha_1}...t_d^{\alpha_d}$.
Notice that $\nu_{a,\rho}(t_it_j)=\nu_{a,\rho}(t_jt_i)=\nu_{a,\rho}(t_i)\nu_{a,\rho}(t_j)$.
Any polynomial $f=\sum_{n\in \Z_+^d}c_nT^n\in K[T_1,...,T_d]_q$
can be written as a finite sum
$f=\sum_{n\in \Z_+^d}c_nt^n+B$ where $\nu_{a,\rho}(B)<\nu_{a,\rho}(\sum_{n\in \Z_+^d}a_nt^n)$.
We see that $\nu_{a,\rho}(f)=max\,|c_n|\rho^n$.
It follows that $\nu_{a,\rho}$
is multiplicative. This concludes the proof. $\blacksquare$

We will say that the seminor $\nu_{a,\rho}$ corresponds to the closed
quantm polydisc $E_q(a,\rho)\subset E_q(0,r)$.
Since the quantum affine space is the union of quantum closed discs
(and hence the Berkovich spectrum of the former is by definition the union
of the Berkovich spectra of the latter) we see that the Berkovich
spectrum of ${\bf A}^d_q$ contains all points $\nu_{a,r}(f)$ with $|a|<r$.

\begin{prp} Previous Proposition holds for the quantum torus
 $T_q(L,\varphi,r)$.

\end{prp}

{\it Proof.}  The proof is basically the same as in the case $q=1$.
Let $f=\sum_{n\in \Z^d}c_nT^n $
be an analytic function
on $T^{an}_q(L,\varphi,r)$. In the notation of the above proof we can rewrite
$f$ as
$f=\sum_{n\in \Z^d}c_n(t+a)^n$, where $(t_i+a_i)^{-1}:=t_i^{-1}\sum_{m\ge 0}(-1)^ma_i^{m}t_i^{-m}$.
Therefore $f=\sum_{m\in \Z_+^d}d_mt^m+\sum_{m\in \Z_-^d}d_mt^m$. Suppose we know that $|d_m|\rho^m\to 0$
as $|m|\to +\infty$. Then similarly to the proof of the previous proposition
one sees that $\nu_{a,\rho}(f):=max_{m}|d_m|\rho^m$ is a multiplicative seminorm. Let us estimate $|d_m|\rho^m$.
For $m\in \Z_+^d$ we have:
$d_m=\sum_{l\in \Z_+^d}b_l^{(m)}c_{l+m}a^{l}$, where $|b_l^{(m)}|\le 1$.
Since $|a|\le \rho$ we have $|d_m|\le max_{l\in \Z_+^d}|c_{l+m}|\rho^{l}\le
max_{l\in \Z_+^d}|c_{l+m}|\rho^{l+m}/\rho^m\le p_m\rho^{-m}$, where $p_m\to 0$ as
$|m|\to \infty$.
Similar estimate holds for $m\in \Z_-^d$.  Therefore $\nu_{a,\rho}(f)$ is a multiplicative
seminorm.
Berkovich spectrum of the quantum torus is, by definition,the union of the Berkovich spectra of
quantum analytic tori $T_q(L,\varphi,r)$ of all multiradii $r$.
Thus we see that
the pair $(a,\rho)$ as above defines a point of
Berkovich spectrum of the quantum analytic torus. $\blacksquare$

\begin{rmk}
It is natural to ask whether
any point of the Berkovich spectrum of the
``commutative" analytic space $(G_m^{an})^d$ appears as a point
of the Berkovich spectrum of the quantum analytic torus, as long  as we choose
$q$ sufficiently close to $1$.
More generally, let us imagine that we have two
quantum affinoid algebras $A$ and $A^{\prime}$ which are addmissible
quotients of the algebras $K\{T_1,...,T_n\}_{Q,r}$
and $K\{T_1,...,T_m\}_{Q^{\prime},r}$ respectively,
where $Q$ and $Q^{\prime}$ are matrices as in Section 3.
One can ask the following question:
for any closed subset $V\subset M(A)$ is there $\varepsilon >0$
such that if $||Q-Q^{\prime}||<\varepsilon$ then
there is a closed subset $V^{\prime}\subset M(A^{\prime})$ homeomorphic to $V$?
Here the norm of the matrix $S=((s_{ij}))$ is defined as
$max_{i,j}\,|s_{ij}|$. If the answer is positive then taking  $Q=id$
we see that the Berkovich spectrum of
an affinoid algebra is a limit of the Berkovich spectrum
of its quantum analytic deformation.

\end{rmk}
The above Proposition shows in a toy-model the drastic difference
with formal deformation quantization. In the latter case
deformations  $A_q$ of a commutative algebra $A_1$
are `` all the same" as long as
$q\ne 1$. In particular, except
of few special values
of $q$, they have the same ``spectra" in the sense of non-commutative
algebraic geometry or representation theory. In analytic case Berkovich spectrum
can contain points which are ``far" from the commutative ones (e.g. we have seen
that the discs $E(a,\rho)\subset E(0,r)$ can be quantized as long as the radius $|\rho|$
is not small). Thus the quantized space contains ``holes" (non-archimedean version
of ``discretization" of the space after quantization).

\section{Non-commutative Stein spaces}

This example is borrowed from [SchT].

Let $K$ be a non-archimedean valuation field, and $A$ be a unital Frechet $K$-algebra.
We say that $A$ is {\it Frechet-Stein} if there is a sequence $v_1\le v_2\le...\le v_n\le...$
of continuous  seminorms on $A$ which define the Frechet topology and:

a) the completions $A_{v_n}$ of the algebras $A/Ker\,v_n$ are left notherian
for all $n\ge 1$;

b) each $A_{v_n}$ is a flat $A_{v_{n+1}}$-module, $n\ge 1$.

Here we do not require that $v_n$ are submultiplicative, only
the inequalities $v_n(xy)\le const\,v_n(x)v_n(y)$. Clearly the sequence
$(A_{v_n})_{n\ge 1}$ is a projective system of algebras and its
projective limit is isomorphic to $A$.

A {\it coherent sheaf} for a Frechet-Stein algebra $(A,v_n)_{n\ge 1}$ is a collection
$M=(M_n)_{n\ge 1}$ such that each $M_n$ is a finite Banach $A_{v_n}$-module,
and for each $n\ge 1$ one has natural isomorphism
$A_{v_n}\otimes_{A_{v_{n+1}}}M_{n+1}\simeq M_n$.

The inverse limit of the projective system $M_n$ is an $A$-module called
the module of {\it global sections} of the coherent sheaf $M$.
Coherent sheaves form an abelian category. It is shown in [SchT] that
the global section functor from coherent sheaves to finite Banach $A$-modules
is exact (this is an analog of $A$ and $B$ theorems of Cartan).

One can define the non-commutative
analytic spectrum of Frechet-Stein algebras in the same way as we did
for Banach algebras,
starting with the category of coherent sheaves.

It is shown in [SchT] that if $G$ is a compact locally analytic group
then  the strong dual $D(G,K)$ to the space of $K$-valued
locally analytic functions is a Frechet-Stein algebra. In the case
$G={\Z}_p$ it is isomorphic to the (commutative) algebra of
power series converging in the open unit disc in the completion
of the algebraic closure of $K$. In general, the Frechet-Stein algebra
structure on $D(G,K)$ is defined by a family of norms which are
submultiplicative only. Coherent sheaves for the algebra $D(G,K)$ 
should give rise
to coherent sheaves on the non-commutative
analytic spectrum $M_{NC}(D(G,K))$, rather than on the Berkovich
spectrum $M(D(G,K))$.

\section{Non-commutative analytic K3 surfaces}

\subsection{General scheme}

The following way of constructing a non-commutative analytic K3 surface $X$ over the field $K={\C}((t))$ was suggested in [KoSo1].

1) We start with a 2-dimensional sphere $B:=S^2$ equipped with an integral affine
structure outside of a finite subset $B^{sing}=\{x_1,...,x_{24}\}$ of $24$ distinct points (see [KoSo1] for the definitions and explanation why $|B^{sing}|=24$).
We assume that the monodromy of the affine structure around each point $x_i$
is conjugate to the $2\times 2$ unipotent Jordan block (it is proved in [KoSo1] that
this restriction enforces the cardinality of $B^{sing}$ to be equal to $24$).

2) In addition to the above data we have an infinite set of ``trees" embedded in $S^2$,
called {\it lines} in [KoSo1]. Precise definition and
the existence of such a set satisfying certain axioms can be found
in [KoSo1], Sections 9, 11.5.

3) The non-commutative analytic space $X_q^{an}$ will be defined defined by a pair
$(X_0,{\cal O}_{X_0,q})$, where $q\in K^{\ast}$ is an arbitrary element satisfying
the condition $|q|=1$, and $X_0$
is a topological K3 surface and ${\cal O}_{X_0,q}$ is a sheaf on
(a certain topology on) $X_0$ of non-commutative noetherian algebras over the field $K$.

4) There is a natural continuous map $\pi:X_0\to S^2$ with the generic fiber
being a two-dimensional torus.
The sheaf ${\cal O}_{X_0,q}$ is uniquely determined by its direct image
${\cal O}_{S^2,q}:=\pi_{\ast}({\cal O}_{X_0,q})$.
Hence the construction of
$X_q$ is reduced to the construction of the sheaf ${\cal O}_{S^2,q}$
on the sphere $S^2$.

5) The sheaf ${\cal O}_{S^2,q}$ is glued from two sheaves: a sheaf
${\cal O}^{sing}_q$ defined in a neighborhood $W$ of the ``singular" subset $B^{sing}$, and a  sheaf ${\cal O}^{nonsing}_q$
on $S^2\setminus W$.

6) The sheaf ${\cal O}^{sing}_q$ is defined by an ``ansatz" described below, while the sheaf
${\cal O}^{nonsing}_q$ is constructed in two steps. First, with any integral affine structure
and an element $q\in K^{\ast}, |q|=1$ one can associate canonically a sheaf ${\cal O}^{can}_q$ of non-commutative algebras over $K$.
Then, with each line $l$ one associates an automorphism $\varphi_l$
of the restriction of ${\cal O}^{can}_q$ to $l$.
The sheaf ${\cal O}^{nonsing}_q$ is obtained
from the restriction of ${\cal O}^{can}_q$ to the complement to the union of all lines by the gluing
procedure by means of $\varphi_l$.

The above scheme was realized in [KoSo1] in the case $q=1$.
In that case one obtains a sheaf ${\cal O}_{X_0,1}$ of
{\it commutative} algebras which is the sheaf of
anayltic functions on the non-archimedean analytic K3 surface. In this section
we will explain what has to be changed in [KoSo1] in order to handle
the case $q\ne 1, |q|=1$.
As we mentioned above, in this case ${\cal O}_{X_0,q}$ will be a sheaf of non-commutative algebras, which is a flat deformation of ${\cal O}_{X_0,0}$. It was observed
in [KoSo1] that ${\cal O}_{X_0,1}$ is a sheaf of Poisson algebras. The sheaf
${\cal O}_{X_0,q}$ is a {\it deformation quantization} of ${\cal O}_{X_0,1}$.
It is an analytic (not a formal) deformation quantization with respect to the
parameter $q-1$. The topology on $X_0$ will be clear from the construction.

\subsection{$\Z$-affine structures and the canonical sheaf}

Let $K$ be a non-archimedean valuation field.
Fix an element
$q\in K^{\ast}, |q|=1$. Let us introduce invertible variables
$\xi, \eta$ such that
$$\eta\xi=q\xi\eta.$$

Then we define a sheaf ${\cal O}^{can}_q$ on $\R^2$ such that for
any open connected subset $U$ one has
$${\cal O}^{can}_q(U)=\left\{
\sum_{n,m\in \Z}c_{n,m}\xi^n \eta ^m\,|\,\forall (x,y)\in U\,\,\,\sup_{n,m}\left(
\log(|c_{n,m}|)+nx+my\right)<\infty\right\}.$$

The above definition is motivated by the following
considerations.

Recall that an integral affine structure ($\Z$-affine structure
for short)
on an $n$-dimensional topological manifold $Y$ is given by
a maximal atlas of  charts such that the change of coordinates
between any two charts is described by the formula
$$ x_i^{\prime}=\sum_{1\le j\le n}a_{ij}x_j+b_i,$$
where $(a_{ij})\in GL(n,{\Z}), (b_i)\in {\R}^n$.
In this case one can speak about the sheaf of $\Z$-affine functions,
i.e. those which can be locally  expressed in affine coordinates by the formula
$f=\sum_{1\le i\le n}a_ix_i+b, a_i\in \Z, b\in \R.$
Another equivalent description: $\Z$-affine structure is given by
a covariant lattice $T^{\Z}\subset TY$ in the tangent bundle (recall that an affine
structure on $Y$ is the same as a torsion free flat connection on the
tangent bundle $TY$).

Let $Y$ be a manifold with $\Z$-affine
structure. The sheaf of $\Z$-affine functions
$Aff_{\Z}:=Aff_{{\Z},Y}$ gives rise to an exact
sequence of sheaves of abelian groups

$$ 0\to {\R}\to Aff_{\Z}\to (T^{\ast})^{\Z}\to 0,$$

where $(T^{\ast})^{\Z}$ is the sheaf associated with the dual to the covariant lattice
$T^{\Z}\subset TY$.

Let us recall the following notion introduced in [KoSo1], Section 7.1.

\begin{defn} A  $K$-affine structure on $Y$
compatible with the given $\Z$-affine structure
is a sheaf $Aff_K$ of abelian groups
on $Y$, an exact sequence of sheaves

$$1\to K^{\times}\to Aff_K\to (T^{\ast})^{\Z}\to 1,$$
together with a homomorphism $\Phi$ of this exact sequence
to the exact sequence of sheaves of abelian groups
$$0\to {\R}\to Aff_{\Z}\to (T^{\ast})^{\Z}\to 0,$$
such that  $\Phi=id$ on $(T^{*})^{\Z}$ and
$\Phi=val$ on $K^{\times}$.

\end{defn}

Since $Y$ carries a $\Z$-affine structure, we
have the corresponding $GL(n,{\Z})\ltimes {\R}^n$-torsor on $Y$,
whose fiber
over a point $x$ consists of all $\Z$-affine coordinate systems
at $x$.

Then one has the following equivalent description of the notion
of $K$-affine structure.

\begin{defn} A $K$-affine structure on $Y$
compatible with the given $\Z$-affine structure is
a $GL(n,{\Z})\ltimes (K^{\times})^n$-torsor
on $Y$ such that the application of $val^{\times n}$ to
$(K^{\times})^n$ gives the initial $GL(n,{\Z})\ltimes {\R}^n$-torsor.
\end{defn}

Assume that $Y$ is oriented and carries a $K$-affine structure compatible
with a given $\Z$-affine structure.
Orientation allows us to reduce to
$SL(n,{\Z})\ltimes (K^{\times})^n$ the structure group of the torsor
defining the $K$-affine structure.
One can define a higher-dimensional version of the sheaf
${\cal O}^{can}_q$ in the following way.
Let $z_1,...,z_n$ be invertible variables such that $z_iz_j=qz_jz_i$,
for all $1\le i<j\le n$. We define the sheaf ${\cal O}^{can}_q$
on ${\R}^n, n\ge 2$ by the same formulas as in the case $n=2$:
$${\cal O}^{can}_q(U)=\left\{
\sum_{I=(I_1,...,I_n)\in {\Z}^n}c_{I}z^I,|\,\forall (x_1,...,x_n)\in U\,\,\,\sup_{I}\left(
\log(|c_{I}|)+\sum_{1\le m\le n}I_mx_m\right)<\infty\right\},$$
where $z^{I}=z_1^{I_1}\dots z_n^{I_n}.$
Since $|q|=1$ the convergency
condition does not depend on the order of variables.

The sheaf ${\cal O}^{can}_q$ can be lifted to $Y$ (we keep the same
notation for the lifting).
In order to do that it suffices to define the action of the group $SL(n,{\Z})\ltimes (K^{\times})^n$
on the canonical sheaf on ${\R}^n$.
Namely, the inverse to an element
$(A,\lambda_1,...,\lambda_n)\in SL(n,{\Z})\ltimes (K^{\times})^n$
acts on monomials as
$$z^{I}=z_1^{I_1}\dots z_n^{I_n}\mapsto \left({\textstyle\prod_{i=1}^n}\lambda_i^{I_i}\right)\,\,z^{A(I)}\,\,\,.$$
The action
of the same element on ${\R}^n$ is given by a similar formula:
$$x=(x_1,\dots ,x_n)\mapsto A(x)-(val(\lambda_1),\dots,val(\lambda_n))\,\,\,.$$

Notice that the stalk of the sheaf ${\cal O}^{can}_q$ over a point
$y\in Y$ is isomorphic to a direct limit of algebras of
functions on quantum analytic tori of various multiradii.

Let $Y=\cup_{\alpha}U_{\alpha}$ be an open covering
by coordinate charts
$U_{\alpha}\simeq V_{\alpha}\subset {\R}^n$
such that for any $\alpha, \beta$ we are given elements
$g_{\alpha,\beta}\in SL(n,{\Z})\ltimes (K^{\times})^n$ satisfying the
$1$-cocycle condition for any triple
$\alpha,\beta,\gamma$. Then the lifting of ${\cal O}^{can}_h$ to $Y$
is obtained via gluing by means of the transformations $g_{\alpha,\beta}$.

Let $G_m^{an}$ be the analytic space corresponding to the
multiplicative group $G_m$ and $(T^n)^{an}:=(G_m^{an})^n$ the $n$-dimensional analytic torus. Then one has a canonically defined
continuous map $\pi_{can}:(T^n)^{an}\to \R^n$ such that
$\pi_{can}(p)=(-val_p(z_1),...,-val_p(z_n))=(log|z_1|_p,...,log|z_n|_p)$,
where $|a|_p$ (resp. $val_p(a)$)   denotes the seminorm (resp. valuation) of an element $a$ corresponding to the point $p$.

For an open  subset $U\in \R^n$ we have a topological $K$-algebra
${\cal O}^{can}_q(U)$ defined by the formulas above.
Notice that a point $x=(x_1,...,x_n)\in \R^n$
defines a multiplicative seminorm $|\sum_Ic_Iz^I|_{exp(x)}$
on the algebra of formal series of $q$-commuting variables
$z_1,...,z_n$ (here $exp(x)=(exp(x_1),...,exp(x_n))$).

Let $M({\cal O}^{can}_q(U))$ be the set of  multiplicative
seminorms $\nu$  on ${\cal O}^{can}_q(U)$ extending
the norm on $K$.
We have defined an embedding $U\to M({\cal O}^{can}_q(U))$,
such that $(x_1,...,x_n)$ corresponds to a seminorm with $|z_i|=exp(x_i),1\le i\le n$.
The map  $\pi_{can}:|\bullet|\mapsto (log|z_1|,...,log|z_n|)$ is a retraction of
$ M({\cal O}^{can}_q(U))$ to the image of $U$.

Let $S_1^n\subset (K^{\times})^n$ be the set of such $(s_1,...,s_n)$ that
$|s_i|=1, 1\le i\le n$. The group $S_1^n$ acts on ${\cal O}^{can}_h(U)$
in such a way that $z_i\mapsto s_iz_i$. Clearly the map $\pi_{can}$ is
$S_1^n$-invariant. For this reason
we will call $\pi_{can}$
a {\it quantum analytic torus fibration} over $U$.
More precisely, we reserve this name for a pair $(U,{\cal O}^{can}_q(U))$,
where the algebra is equipped with the $S_1^n$-action. We suggest to think
about such a pair as of the algebra ${\cal O}_q(\pi^{-1}_{can}(U))$ of analytic
functions on the open subset $\pi^{-1}_{can}(U)$ of
the non-commutative analytic torus $(T^n_q(\Z^n)^{an},\varphi_0)$,
where $\varphi_0((a_1,...,a_n),(b_1,...,b_n))=a_1b_1+...+a_nb_n$.

We can make a category of the above pairs, defining a morphism
$(U,{\cal O}_q^{can}(U))\to (V,{\cal O}_q^{can}(V))$ as a pair
$(f,\phi)$ where $f:U\to V$ is a continuous map and
$\phi:{\cal O}_q^{can}(f^{-1}(V))\to {\cal O}_q^{can}(V)$ is a
$S_1^n$-equivariant homomorphism of algebras.
In particular we have the notion of isomorphism of quantum
analytic torus fibrations.

Let $U\subset \R^n$ be an open set and $A$ be a non-commutative
affinoid $K$-algebra equipped with
a $S_1^n$-action. We say that a pair $(U,A)$ defines a quantum analytic torus
fibration over $U$ if it is isomorphic to the pair
$(U,{\cal O}_q^{can}(U))$. Notice that morphisms of quantum analytic torus fibrations
are compatible with the restrictions on the open subsets. Therefore we can introduce
a topology on $(T^n)^{an}$ taking $\pi_{can}^{-1}(U), U\subset \R^2$ as open subsets,
and make a ringed space assigning the algebra ${\cal O}_q(U)$ to the open set
$\pi_{can}^{-1}(U)$. We will denote the sheaf by $(\pi_{can})_{\ast}({\cal O}_{(T^n)^{an},q})$.
Its global sections (for $U=\R^n$) coincides with the projective limit of algebras
of analytic functions
on  quantum tori of all possible multiradii. Slightly abusing the terminology we will call
the above ringed space a quantum analytic torus.

\begin{rmk} The above definition is a toy-model of a
$q$-deformation of the natural retraction $X^{an}\to Sk(X)$ of
the analytic space $X^{an}$  associated with the maximally degenerate Calabi-Yau manifold $X$
onto its skeleton $Sk(X)$ (see [KoSo1]). Analytic torus fibrations introduced in [KoSo1]
are ``rigid analytic" analogs of Lagrangian torus fibrations in symplectic geometry.
Moreover, the mirror symmetry functor (or rather its incarnation as a Fourier-Mukai
transform) interchanges these two types of torus fibrations for mirror dual
Calabi-Yau manifolds.

\end{rmk}

It turns out that the  sheaf ${\cal O}_q^{can}$ is not good for construction of
a non-commutative analytic K3 surface. We will explain later
how it should be modified. Main reason for the complicated modification
procedure comes from Homological Mirror Symmetry, as
explained in [KoSo1]. In a few words, the derived category of
coherent sheaves on a non-commutative analytic K3 surface should
be equivalent to a certain deformation of the Fukaya category
of the mirror dual K3 surface. If the K3 surface is realized
as an elliptic fibration over ${\C}{\bf P}^1$ then there are
fibers (they are $2$-dimensional Lagrangian
tori) which contain boundaries
of holomorphic discs. Those discs give rise to an infinite set
of lines on the base of the fibration. In order to have the above-mentioned
categorical equivalence one should modify the canonical sheaf for
each line.

\subsection{Model near a singular point}

Let us fix $q\in K^{\ast}, |q|=1$.

We start with the open covering of $\R^2$ by the following sets $U_i, 1\le i\le 3$.
Let us fix a number $0<\varepsilon<1$ and define
$$\begin{array} {lll}
U_1 & = & \{(x,y)\in {\R}^2|x<\varepsilon |y|\,\}\\
U_2 & = & \{(x,y)\in {\R}^2|x>0, y<\varepsilon x\,\}\\
U_3 & = & \{(x,y)\in {\R}^2|x>0, y>0\}
\end{array}$$
Clearly ${\R}^2\setminus\{(0,0)\}=U_1\cup U_2\cup U_3$.
We will also need a slightly modified domain $U_2'\subset U_2$ defined as
$\{(x,y)\in {\R}^2|x>0, y<\frac{\varepsilon}{1+\varepsilon} x\,\}$.

Recall that one has a canonical map $\pi_{can}:(T^2_q)^{an}\to {\R}^2.$

We define $T_i:=\pi_{can}^{-1}(U_i), i=1,3$ and
 $T_2:=\pi_{can}^{-1}(U_2')$ (see the explanation below).
Then the projections $\pi_i: T_i\to U_i$ are given by the formulas
$$
\pi_i(\xi_i,\eta_i)=\pi_{can}(\xi_i,\eta_i)=
(\log|\xi_i|, \log|\eta_i|), \,\,\,i=1,3$$
$$\pi_2(\xi_2,\eta_2)=\left\{\begin{array}{ll}(\log|\xi_2|, \log|\eta_2|)& \mbox{  if }|\eta_2|<1\\
(\log|\xi_2|-\log|\eta_2|, \log|\eta_2|) & \mbox{ if }
|\eta_2|\ge 1
\end{array}\right. $$
 In these formulas $(\xi_i,\eta_i)$ are  coordinates
on $T_i, 1\le i\le 3$. More pedantically, one should say that
for each $T_i$ we are given an algebra ${\cal O}_q(T_i)$ of series
$\sum_{m,n}c_{mn}\xi_i^m\eta_i^n$ such that $\xi_i\eta_i=q\eta_i\xi_i$,
and for a seminorm $|\bullet|$ corresponding to a point of $T_i$ (which means
that $(log|\xi_i|,log|\eta_i|)\in U_i$) one has:
$sup_{m,n}(m\,log|\xi_i|+n\,log|\eta_i|)<+\infty$. In this way we obtain a sheaf
of non-commutative algebras on the set $U_i$, which is the subset of the  Berkovich spectrum
of the algebra
${\cal O}_q(T_i)$. We will denote this sheaf by ${\cal O}_{T_i,q}$.

Let us introduce the sheaf ${\cal O}^{can}_q$ on ${\R}^2\setminus \{(0,0)\}$.
It is defined as $(\pi_i)_*\left(\O_{T_i,q}\right)$ on each domain $U_i$,
with identifications
$$\begin{array}{llcl}
(\xi_1,\eta_1) & = & (\xi_2,\eta_2) & \mbox{ on } U_1\cap U_2 \\
(\xi_1,\eta_1) & = & (\xi_3,\eta_3) & \mbox{ on } U_1\cap U_3 \\
(\xi_2,\eta_2) & = & (\xi_3\eta_3,\eta_3) & \mbox{ on } U_2\cap U_3
\end{array}$$

Let us introduce the sheaf ${\cal O}_q^{sing}$ on ${\R}^2\setminus \{(0,0)\}$.
On the sets $U_1$ and $U_2\cup U_3$ this
sheaf is isomorphic to  ${\cal O}^{can}_q$
(by identifying of coordinates $(\xi_1,\eta_1)$
and of glued coordinates $(\xi_2,\eta_2)$ and $(\xi_3,\eta_3)$ respectively).
On the intersection $U_1\cap (U_2\cup U_3)$ we identify two
copies of  the canonical sheaf by an automorphism $\varphi$ of
${\cal O}^{can}_q$.
More precisely, the automorphism is given (we skip the index of the
coordinates) by
 $$\varphi(\xi,\eta)=\left\{\begin{array}{cll}
(\xi(1+\eta),\eta) & \mbox{ on } & U_1\cap U_2 \\
(\xi(1+1/\eta),\eta) & \mbox{ on } & U_1\cap U_3
\end{array}
\right.$$

\subsection{Lines and automorphisms}

We refer the reader to [KoSo1] for the precise definition of
the set of lines and axioms
this set is required to obey. Roughly speaking, for a manifold $Y$
which carries a $\Z$-affine structure a line $l$ is defined by
a continuous map $f_l: (0,+\infty)\to Y$ and a covariantly constant
nowhere vanishing integer-valued $1$-form
$\alpha_l\in \Gamma((0,+\infty),f_l^{\ast} ((T^\ast)^\Z)$.
A set ${\cal L}$ of lines is required to be decomposed into a disjoint union
${\cal L}={\cal L}_{in}\cup {\cal L}_{com}$ of {\it initial}
and {\it composite} lines. Each composite line is obtained as a result
of a finite number of ``collisions" of initial lines. A collision
is described by a $Y$-shape figure, where the leg of $Y$ is a composite line,
while two other seqments are ``parents" of the leg. A construction of the
set ${\cal L}$ satisfying the axioms from [KoSo1] was proposed in
[KoSo1], Section 9.3.

With each line $l$ we can associate a continuous family of automorphisms
of stalks of sheaves of algebras
$\varphi_l(t): ({\cal O}^{can}_q)_{Y,f_l(t)}\to ({\cal O}^{can}_q)_{Y,f_l(t)}$.

Automorphisms $\varphi_l$ can be defined in the following way (see [KoSo1], Section 10.4).

First we choose affine coordinates in a neighborhhod of a point
$b\in B\setminus B^{sing}$, identifyin $b$ with the
point $(0,0)\in \R^2$.
Let $l=l_+\in {\cal L}_{in}$ be (in the standard affine coordinates)
a line in the half-plane $y>0$ emerging from $(0,0)$
(there is another such line $l_{-}$
in the half-plane $y<0$, see [KoSo1] for the details).
Assume that $t$
is sufficiently small. Then we define
$\varphi_l(t)$
on topological generators $\xi,\eta$ by the formula

$$\varphi_l(t)(\xi,\eta)=(\xi(1+1/\eta),\eta).$$

In order to extend $\varphi_l(t)$ to the interval $(0,t_0)$, where
$t_0$ is not small, we cover the corresponding segment of $l$ by
open charts. Then a change of affine coordinates transforms
$\eta$ into a monomial multiplied by a constant from $K^{\times}$.
Moreover, one can choose the change of coordinates in such a way that
$\eta\mapsto C\eta$ where $C\in K^{\times}, |C|<1$ (such change
of coordinates preserve the $1$-form $dy$. Constant $C$ is equal to
$exp(-L)$, where $L$ is the length of the segment of $l$ between
two points in different coordinate charts).
Therefore $\eta$ extends analytically in a unique way to an element of
$\Gamma((0,+\infty), f_l^{\ast}(({\cal O}^{can}_q)^{\times}))$.
Moreover the norm $|\eta|$ strictly decreases as $t$ increases,
and remains strictly smaller than $1$. Similarly to [KoSo1], Section 10.4
one deduces that $\varphi_l(t)$ can be extended for all $t>0$.
This defines $\varphi_l(t)$ for $l\in {\cal L}_{in}$.

Next step is to extend $\varphi_l(t)$ to the case when $l\in {\cal L}_{com}$,
i.e. to the case when the line is obtained as a result of a collision
of two lines belonging to ${\cal L}_{in}$.
Following [KoSo1], Section 10, we introduce a group $G$ which contains
all the automorphisms $\varphi_l(t)$, and then prove the factorization theorem
(see [KoSo1], Theorem 6) which allows us to define $\varphi_l(0)$ in the case
when $l$ is obtained as a result of a collision of two lines $l_1$ and $l_2$.
Then we extend $\varphi_l(t)$ analytically for all $t>0$ similarly to the
case $l\in {\cal L}_{in}$.

More precisely, the construction of $G$ goes such as follows.
Let $(x_0,y_0)\in \R^2$ be a point,
$\alpha_1,\alpha_2\in (\Z^2)^\ast$ be $1$-covectors such that
$\alpha_1\wedge\alpha_2>0$.
Denote by $V=V_{(x_0,y_0),\alpha_1,\alpha_2}$ the closed angle
$$\{(x,y)\in \R^2|\langle \alpha_i, (x,y)-(x_0,y_0)\rangle\ge 0, i=1,2\,\}$$

Let ${\cal O}_q(V)$ be a $K$-algebra consisting
of series $f=\sum_{n,m \in \Z}c_{n,m}\xi^{n}\eta^{m}$,
such that $\xi\eta=q\eta\xi$ and $c_{n,m}\in K$ satisfy the condition that
for all $(x,y)\in V$ we have:

\begin{enumerate}
\item if $c_{n,m}\ne 0$ then $\langle (n,m),(x,y)-(x_0,y_0)\rangle\le 0$,
where we identified $(n,m)\in \Z^2$
with a covector in $(T_{p}^{\ast}Y)^{\Z}$;
\item $\log|c_{n,m}|+nx+my\to -\infty$ as long as $|n|+|m|\to +\infty$.
\end{enumerate}

For an integer covector $\mu=adx+bdy\in (\Z^2)^*$ we denote
by $R_{\mu}$ the monomial $\xi^a\eta^b$.
Then we consider a pro-unipotent group $G:=G(q,\alpha_1,\alpha_2,V)$ of automorphisms of
${\cal O}_q(V)$  having the form

$$g=
\sum_{n_1,n_2\ge 0,n_1+n_2>0}c_{n_1,n_2}R_{\alpha_1}^{-n_1}R_{\alpha_2}^{-n_2} $$
where
$$\log|c_{n,m}|-n_1\langle \alpha_1, (x,y)\rangle-n_2\langle \alpha_2, (x,y)\rangle\le 0\,\,\,\,\forall \,(x,y)\in V$$
The latter condition is equivalent to
$\log|c_{n,m}|-\langle n_1 \alpha_1+n_2\alpha_2, (x_0,y_0)\rangle\le 0$.

Fixing the ratio $\lambda=n_2/n_1\in [0,+\infty]_{\Q}:=\Q_{\ge 0}\cup\infty$ we obtain a subgroup $G_{\lambda}:=G_{\lambda}(q,\alpha_1,\alpha_2,V)\subset G$.
There is a natural map $\prod_{\lambda}G_{\lambda}\to G$, defined as in [KoSo1],
Section 10.2. The above-mentioned factorization theorem states that this
map is a bijection of sets.

Let us now assume that lines $l_1$ and $l_2$ collide at
$p=f_{l_1}(t_1)=f_{l_2}(t_2)$,
generating the line $l\in {\cal L}_{com}$.
Then $\varphi_l(0)$ is defined with the help of factorization theorem.
More precisely, we set $\alpha_i:=\alpha_{l_i}(t_i),\,\,i=1,2$ and the angle
$V$ is the intersection of certain half-planes $P_{l_1,t_1}\cap P_{l_2,t_2}$
defined in [KoSo1], Section 10.3. The half-plane $P_{l,t}$ is contained
in the region of convergence of $\varphi_l(t)$.
By construction, the elements
$g_0:=\varphi_{l_1}(t_1)$  and $g_{+\infty}:=\varphi_{l_2}(t_2)$ belong respectively
to $G_0$ and $G_{+\infty}$.
The we have:
$$g_{+\infty} g_0= {\textstyle \prod_\to}\left(
(g_\lambda)_{\lambda\in [0,+\infty]_\Q}\right)=
g_0\dots g_{1/2}\dots g_1 \dots
g_{+\infty}.$$

Each term $g_\lambda$ with $0<\lambda= n_1/n_2 <+\infty$
corresponds to the newborn line $l$ with the
direction covector
$n_1\alpha_{l_1}(t_1)+n_2\alpha_{l_2}(t_2)$.
Then we set  $\varphi_l(0):=g_{\lambda}$.
This transformation is defined by a series
which is convergent in a neighborhood of $p$, and using the analytic continuation
we obtain $\varphi_l(t)$ for $t>0$, as we said above. Recall that
every line carries an integer $1$-form $\alpha_l=adx+bdy$. By construction,
$\varphi_l(t)\in G_{\lambda}$, where $\lambda$ is the slope of $\alpha_l$.

Having automorphisms $\varphi_l$ assigned to lines $l\in {\cal L}$ we
proceed as in [KoSo1], Section 11, modifying the sheaf
${\cal O}^{can}_q$ along each line. We denote the resulting
sheaf ${\cal O}^{nonsing}_q$. By construction it is isomorphic to the sheaf
${\cal O}^{sing}_q$ in a neighborhood of the point $(0,0)$.

Let us now consider the manifold $Y=S^2\setminus B^{sing}$, i.e.
the complement of $24$ points on the sphere $S^2$ equipped with
the $\Z$-affine structure, which has standard singularity at
each point $x_i\in B^{sing}, 1\le i\le 24$ (see Section 7.1).
Using the above construction (with any choice of set of lines
on $S^2$) we define the sheaf ${\cal O}_{S^2,q}^{nonsing}$ on $Y$.
Notice that in a small neighborhood of each singular point $x_i$
the sheaf ${\cal O}_{S^2,q}^{nonsing}$ is isomorphic to the sheaf
${\cal O}_q^{sing}$ (in fact they become isomorphic
after identification of the punctured neigborhood of $x_i$ with
the punctured neighborhood of $(0,0)\in \R^2$ equipped with the
standard singular $\Z$-affine structure (see Section 7.1 and [KoSo1],
Section 6.4 for the description of the latter).
In the next subsection we will give an alternative description of
the sheaf ${\cal O}_q^{sing}$. It follows from that description
that ${\cal O}_q^{sing}$ can be extended to the point $(0,0)$.
It gives a sheaf ${\cal O}_{S^2,q}^{sing}$ in the neighborhood of $B^{sing}$.
As a result we will obtain the sheaf ${\cal O}_{S^2,q}$ of
non-commutative $K$-algebras on the whole sphere $S^2$
such that it is isomorphic to ${\cal O}_{S^2,q}^{nonsing}$ on the
complement of $B^{sing}$ and isomorphic to ${\cal O}_{S^2,q}^{sing}$
in a neighborhood of $B^{sing}$.

\subsection{About the sheaf ${\cal O}^{sing}_q$}

We need to check that the sheaf ${\cal O}_{S^2,q}$ is a flat deformation of
the sheaf ${\cal O}_{S^2}$ constructed in [KoSo1].
For the sheaf ${\cal O}^{nonsing}_{S^2,q}$ this follows from the construction.
Indeed, the algebra of analytic functions on the quantum analytic torus (of any multiradius)
is a flat deformation of the algebra of analytic functions on the corresponding
``commutative" torus, equipped with the Poisson bracket $\{x,y\}=xy$.
The group $G=G(q)$ described in the previous subsection is a flat deformation
of its ``commutative" limit $G(1):=G(q=1)$ defined in [KoSo1], Section 10.
The group $G(1)$ preserves the above Poisson bracket.

In order to complete the construction of the non-commutative space
analytic K3 surface $X_h$
we need to investigate the sheaf ${\cal O}_q^{sing}$ and prove that it is
a flat deformation with respect to $q-1$ of the sheaf ${\cal O}^{model}$
introduced in [KoSo1], Section 8. First we recall the definition
of the latter.

Let $S\subset {\bf A}^3$ be an algebraic surface given by equation
$ (\alpha\beta-1)\gamma=1$ in coordinates $(\alpha,\beta,\gamma)$,
and $S^{an}$ be the corresponding analytic space.
We define a continuous map
$f: S^{an}\to {\R}^3$ by the formula
$f(\alpha,\beta,\gamma)=(a,b,c)$
where $a=\max(0,\log|\alpha|_p), b=\max (0, \log|\beta|_p),
c= \log|\gamma|_p=-\log|\alpha\beta-1|_p$. Here
$|\cdot|_p=\exp(-val_p(\cdot))$ denotes the mulitplicative seminorm
corresponding to the point $p\in S^{an}$.

Let us consider the
embedding $j: {\R}^2\to {\R}^3$ given by  formula
$$j(x,y)=\left\{\begin{array}{lll}
(-x\,,\, \max(x+y,0)\,,\, -y\,) &
\mbox{ if } & x\le 0\\
 (\,0\,,\, x+\max(y,0)\,,\, -y\,) & \mbox{ if } & x\ge 0
\end{array}\right.$$
One can easily check that the image of $j$ coincides with
the image of $f$.
Let us denote by $\pi: S^{an}\to {\R}^2$ the map $j^{(-1)}\circ f$.
Finally, we denote $\pi_{\ast}({\cal O}_{S^{an}})$ by
${\cal O}^{model}:={\cal O}^{model}_{\R^2}$.
It was shown in [KoSo1], Section 8, that ${\cal O}^{model}$
is canonically isomorphic to the sheaf ${\cal O}_{q=1}^{sing}$
(the latter is defined as a modification of the sheaf ${\cal O}^{can}_{q=1}$
by means of the automorphism ${\varphi}$, given by the formula at the end
of Section 7.3 for commuting variables $\xi$ and $\eta$).

Let us consider a non-commutative $K$-algebra $A_q(S)$ generated
by generators $\alpha,\beta,\gamma$ subject to the following
relations:
$$\alpha\gamma=q\gamma\alpha,\beta\gamma=q\gamma\beta,$$
$$\beta\alpha-q\alpha\beta=1-q,$$
$$(\alpha\beta-1)\gamma=1.$$
For $q=1$ this algebra concides with the algebra of regular
functions on the surface $X\subset {\bf A}^3_K$ given
by the equation $(\alpha\beta-1)\gamma=1$ and moreover, it is a flat
deformation of the latter with respect to the parameter $q-1$.

Recall that in Section 7.3 we defined three open subsets $T_i, 1\le i\le 3$
of the two-dimensional quantum analytic torus $(T^2)^{an}$, The subset $T_i$
is defined as a ringed space $(\pi_i^{-1}(U_i), {\cal O}_{T_i,q})$, where
$U_i$ are open subsets of $\R^2$ and ${\cal O}_{T_i,q}$ is a sheaf of
non-commutative algebras, uniquely determined by the $K$-algebra
${\cal O}_q(T_i)$ of its global sections.

We define morphisms $g_i: T_i\mono S, 1\le i\le 3$ by the following formulas

$$\begin{array}{lll}
g_1(\xi_1,\eta_1) & = & ({{\xi_1}^{-1}}, \xi_1(1+\eta_1), {\eta_1}^{-1})\\
g_2(\xi_2,\eta_2) & = & ({(1+\eta_2){\xi_2}^{-1}}, \xi_2, {\eta_2}^{-1})\\
g_3(\xi_3,\eta_3) & = &
((1+\eta_3)(\xi_3\eta_3)^{-1}, \xi_3\eta_3, (\eta_3)^{-1})
\end{array}$$

Pedantically speaking this means that for each $1\le i\le 3$ we have a homomorphim
of $K$-algebras $A_q(S)\to {\cal O}_q(T_i)$ such that $\alpha$
is mapped to the first coordinate of $g_i$, $\beta$ is mapped to the second
coordinate and $\gamma$ is mapped to the third coordinate.
One checks directly that three coordinates obey the relations
between $\alpha,\beta,\gamma$.
Modulo $(q-1)$ these morphisms are inclusions.
In the non-commutative case they induce embeddings of $M({\cal O}_q(T_i)), 1\le i\le 3$
into the set $M(A_q(S))$ of multiplicative seminorms on $A_q(S)$.

Notice that in the commutative case we have:
$j\circ \pi_i=f\circ g_i$ and
$f^{-1}(j(U_i))=g_i(T_i)$ for all $1\le i\le 3$.
Using this observation we can decompose a neighborhood $V$
of $\pi^{-1}(0,0)$ in $S^{an}$ into three open analytic
subspaces and describe explicitly algebras of analytic functions as series
in coordinates $(\alpha,\beta)$ or $(\beta,\gamma)$ or $(\alpha,\gamma)$
(choice of the coordinates depend on the domain) with certain grows condtions
on the coefficients of the series. This gives explicit description
of the algebra $\pi_{\ast}({\cal O}_{S^{an}}(\pi(V))$. Then we
{\it declare} the same description in the non-commutative case to be
the answer for the direct image. Non-commutativity does not affect
the convergency condition because $|q|=1$. This description, perhaps,
can be obtained from the ``general theory" which will developed elsewhere.
The direct check, as in the commutative case, shows the compatibility of this
description of the direct image sheaf with the description of ${\cal O}^{nonsing}_q$ in the neighborhood
of $(0,0)$. Therefore we can glue both sheaves together, obtaining
${\cal O}_{S^2,q}$. This concludes the construction.

\vspace{5mm}

{\bf Bibliography}

\vspace{3mm}

[BGR] S. Bosch, U. G\"untzer, R. Remmert, Non-archimedean analysis. Springer-Verlag, 1984.

\vspace{2mm}

[Be1] V. Berkovich, Spectral theory and analytic geometry over non-archimedean fields.
AMS Math. Surveys and Monographs, n. 33, 1990.

\vspace{2mm}

[Be2] V. Berkovich, Etale cohomology for non-Archimedean analytic spaces,
Publ. Math. IHES 78 (1993), 5-161.

\vspace{2mm}

[Be3] V. Berkovich, Smooth p-adic analytic spaces are locally contractible, Inv. Math., 137 (1999),
1-84.
\vspace{2mm}

[Co] A. Connes, Non-commutative geometry, Academic Press, 1994.

\vspace{2mm}

[FG] V. Fock, A. Goncharov, Cluster ensembles, quantization and
the dilogarithm, math.AG/0311245.

\vspace{2mm}

[FvP] J. Fresnel, M. van der Put, Rigid analytic geometry and applications. Birkhauser, 2003.

\vspace{2mm}

[KoSo1] M. Kontsevich, Y. Soibelman, Affine structures and non-archimedean analytic spaces,
math.AG/0406564.

\vspace{2mm}

[KoSo2] M. Kontsevich, Y. Soibelman,Homological mirror symmetry and torus fibrations,
math.SG/0011041.
\vspace{2mm}

[KoR] M. Kontsevich, A. Rosenberg, Non-commutative smooth spaces, math.AG/9812158.

\vspace{2mm}

[KorSo] L. Korogodsky, Y. Soibelman, Alegbras of functions on quantum groups. I, Amer. Math. Soc., 1997.

\vspace{2mm}

[Lap] Lapchik, personal communications, 2004-2005, see also www.allalapa.com

\vspace{2mm}

[R1] A. Rosenberg, Non-commutative algebraic geometry and representations of quantized algebras.
Kluwer Academic Publishers, 1995.

\vspace{2mm}

[RSo] A. Rosenberg, Y. Soibelman, Non-commutative analytic spaces, in preparation.

\vspace{2mm}

[SchT] P. Schneider, J. Teitelbaum, Algebras of p-adic distributions and admissible representations,
Invent. math. 153, 145 - 196 (2003)
\vspace{2mm}

[SoVo] Y. Soibelman, V. Vologodsky, Non-commutative compactifications and elliptic curves,
math.AG/0205117, published in Int. Math.Res. Notes, 28 (2003).

\vspace{5mm}

Y.S.: {email: soibel@math.ksu.edu}

\vspace{2mm}

Address: Department of Mathematics, Kansas State University, Manhattan, KS 66506, USA

\end{document}